\documentclass[11pt]{article}
\usepackage{textcomp}
\usepackage{amsbsy}
\usepackage{latexsym}
\usepackage[mathscr]{eucal}
\usepackage{amsfonts,amsmath,amsthm}
\usepackage{amssymb}
\usepackage[usenames]{xcolor}

\usepackage{fullpage}

\begin{document}
\bibliographystyle{plain}
\title
{
Diffusion Coefficients Estimation for Elliptic Partial Differential Equations
}
\author{ 
Andrea Bonito, Albert Cohen, Ronald DeVore, Guergana Petrova, and Gerrit Welper
\thanks{%
    This research was supported by the ONR Contracts  N00014-15-1-2181 and N0014-16-2706; 
  the  NSF Grants
       DMS 1521067 and DMS 1254618; the DARPA Grant  HR0011619523 through Oak Ridge National Laboratory,
    }  }

\newcommand{\lb}[1]{$\mathbf{PC}(#1)$}
\newcommand{\grad}{$\mathbf{GC}(\theta,M)$}

\hbadness=10000
\vbadness=10000
\newtheorem{lemma}{Lemma}[section]
\newtheorem{proposition}[lemma]{Proposition}
\newtheorem{corollary}[lemma]{Corollary}
\newtheorem{theorem}[lemma]{Theorem}
\newtheorem{remark}{Remark}[section]
\newtheorem{example}[lemma]{Example}
\newtheorem{definition}[lemma]{Definition}
\newtheorem{assumption}[lemma]{Assumption}
\newtheorem{result}{Result}

%
\newenvironment{disarray}{\everymath{\displaystyle\everymath{}}\array}{\endarray}

\def\RR{\rm \hbox{I\kern-.2em\hbox{R}}}
\def\NN{\rm \hbox{I\kern-.2em\hbox{N}}}
\def\ZZ{\rm {{\rm Z}\kern-.28em{\rm Z}}}
\def\CC{\rm \hbox{C\kern -.5em {\raise .32ex \hbox{$\scriptscriptstyle
|$}}\kern
-.22em{\raise .6ex \hbox{$\scriptscriptstyle |$}}\kern .4em}}
\def\vp{\varphi}
\def\<{\langle}
\def\>{\rangle}
\def\t{\tilde}
\def\i{\infty}
\def\e{\varepsilon}
\def\sm{\setminus}
\def\nl{\newline}
\def\o{\overline}
\def\wt{\widetilde}
\def\wh{\widehat}
\def\cT{{\cal T}}
\def\cA{{\cal A}}
\def\cI{{\cal I}}
\def\cV{{\cal V}}
\def\cB{{\cal B}}
\def\cF{{\cal F}}
\def\cY{{\cal Y}}

\def\cD{{\cal D}}
\def\cP{{\cal P}}
\def\cJ{{\cal J}}
\def\cM{{\cal M}}
\def\cO{{\cal O}}
\def\Chi{\raise .3ex
\hbox{\large $\chi$}} \def\vp{\varphi}
\def\lsima{\hbox{\kern -.6em\raisebox{-1ex}{$~\stackrel{\textstyle<}{\sim}~$}}\kern -.4em}
\def\lsim{\hbox{\kern -.2em\raisebox{-1ex}{$~\stackrel{\textstyle<}{\sim}~$}}\kern -.2em}
\def\[{\Bigl [}
\def\]{\Bigr ]}
\def\({\Bigl (}
\def\){\Bigr )}
\def\[{\Bigl [}
\def\]{\Bigr ]}
\def\({\Bigl (}
\def\){\Bigr )}
\def\L{\pounds}
\def\pr{{\rm Prob}}
\newcommand{\cs}[1]{{\color{magenta}{#1}}}
\def\ds{\displaystyle}
\def\ev#1{\vec{#1}}     
\newcommand{\lt}{\ell^{2}(\nabla)}
\def\Supp#1{{\rm supp\,}{#1}}
\def\R{\mathbb{R}}
\def\E{\mathbb{E}}
\def\nl{\newline}
\def\T{{\relax\ifmmode I\!\!\hspace{-1pt}T\else$I\!\!\hspace{-1pt}T$\fi}}
\def\N{\mathbb{N}}
\def\Z{\mathbb{Z}}
\def\N{\mathbb{N}}
\def\Zd{\Z^d}
\def\Q{\mathbb{Q}}
\def\C{\mathbb{C}}
\def\Rd{\R^d}
\def\gsim{\mathrel{\raisebox{-4pt}{$\stackrel{\textstyle>}{\sim}$}}}
\def\sime{\raisebox{0ex}{$~\stackrel{\textstyle\sim}{=}~$}}
\def\lsim{\raisebox{-1ex}{$~\stackrel{\textstyle<}{\sim}~$}}
\def\div{\mbox{div}}
\def\M{M}  \def\NN{N}                  
\def\L{{\ell}}               
\def\Le{{\ell^1}}            
\def\Lz{{\ell^2}}
\def\Let{{\tilde\ell^1}}     
\def\Lzt{{\tilde\ell^2}}
\def\Ltw{\ell^\tau^w(\nabla)}
\def\t#1{\tilde{#1}}
\def\la{\lambda}
\def\La{\Lambda}
\def\ga{\gamma}
\def\BV{{\rm BV}}
\def\Ga{\eta}
\def\al{\alpha}
\def\cZ{{\cal Z}}
\def\cA{{\cal A}}
\def\cU{{\cal U}}
\def\argmin{\mathop{\rm argmin}}
\def\argmax{\mathop{\rm argmax}}
\def\prob{\mathop{\rm prob}}

\def\cO{{\cal O}}
\def\cA{{\cal A}}
\def\cC{{\cal C}}
\def\cS{{\cal F}}
\def\bu{{\bf u}}
\def\bz{{\bf z}}
\def\bZ{{\bf Z}}
\def\bI{{\bf I}}
\def\cE{{\cal E}}
\def\cD{{\cal D}}
\def\cG{{\cal G}}
\def\cI{{\cal I}}
\def\cJ{{\cal J}}
\def\cM{{\cal M}}
\def\cN{{\cal N}}
\def\cT{{\cal T}}
\def\cU{{\cal U}}
\def\cV{{\cal V}}
\def\cW{{\cal W}}
\def\cL{{\cal L}}
\def\cB{{\cal B}}
\def\cG{{\cal G}}
\def\cK{{\cal K}}
\def\cX{{\cal X}}
\def\cS{{\cal S}}
\def\cP{{\cal P}}
\def\cQ{{\cal Q}}
\def\cR{{\cal R}}
\def\cU{{\cal U}}
\def\bL{{\bf L}}
\def\bl{{\bf l}}
\def\bK{{\bf K}}
\def\bC{{\bf C}}
\def\X{X\in\{L,R\}}
\def\ph{{\varphi}}
\def\D{{\Delta}}
\def\H{{\cal H}}
\def\bM{{\bf M}}
\def\bx{{\bf x}}
\def\bj{{\bf j}}
\def\bG{{\bf G}}
\def\bP{{\bf P}}
\def\bW{{\bf W}}
\def\bT{{\bf T}}
\def\bV{{\bf V}}
\def\bv{{\bf v}}
\def\bt{{\bf t}}
\def\bz{{\bf z}}
\def\bw{{\bf w}}
\def \span{{\rm span}}
\def \meas {{\rm meas}}
\def\rhom{{\rho^m}}
\def\diff{\hbox{\tiny $\Delta$}}
\def\EE{{\rm Exp}}
\def\lll{\langle}
\def\argmin{\mathop{\rm argmin}}
\def\codim{\mathop{\rm codim}}
\def\rank{\mathop{\rm rank}}

\def\argmax{\mathop{\rm argmax}}
\def\dJ{\nabla}
\newcommand{\ba}{{\bf a}}
\newcommand{\bb}{{\bf b}}
\newcommand{\bc}{{\bf c}}
\newcommand{\bd}{{\bf d}}
\newcommand{\bs}{{\bf s}}
\newcommand{\bff}{{\bf f}}
\newcommand{\bp}{{\bf p}}
\newcommand{\bg}{{\bf g}}
\newcommand{\by}{{\bf y}}
\newcommand{\br}{{\bf r}}
\newcommand{\be}{\begin{equation}}
\newcommand{\ee}{\end{equation}}
\newcommand{\bea}{$$ \begin{array}{lll}}
\newcommand{\eea}{\end{array} $$}
\def \Vol{\mathop{\rm  Vol}}
\def \mes{\mathop{\rm mes}}
\def \Prob{\mathop{\rm  Prob}}
\def \exp{\mathop{\rm    exp}}
\def \sign{\mathop{\rm   sign}}
\def \sp{\mathop{\rm   span}}
\def \rad{\mathop{\rm   rad}}
\def \vphi{{\varphi}}
\def \csp{\overline \mathop{\rm   span}}
%
%
\newcommand{\beqn}{\begin{equation}}
\newcommand{\eeqn}{\end{equation}}
\def\beginproof{\noindent{\bf Proof:}~ }
\def\endproof{\hfill\rule{1.5mm}{1.5mm}\\[2mm]}

\newenvironment{Proof}{\noindent{\bf Proof:}\quad}{\endproof}

\renewcommand{\theequation}{\thesection.\arabic{equation}}
\renewcommand{\thefigure}{\thesection.\arabic{figure}}

\makeatletter
\@addtoreset{equation}{section}
\makeatother

\newcommand\abs[1]{\left|#1\right|}
\newcommand\clos{\mathop{\rm clos}\nolimits}
\newcommand\trunc{\mathop{\rm trunc}\nolimits}
\renewcommand\d{d}
\newcommand\dd{d}
\newcommand\diag{\mathop{\rm diag}}
\newcommand\dist{\mathop{\rm dist}}
\newcommand\diam{\mathop{\rm diam}}
\newcommand\cond{\mathop{\rm cond}\nolimits}
\newcommand\eref[1]{{\rm (\ref{#1})}}
\newcommand{\iref}[1]{{\rm (\ref{#1})}}
\newcommand\Hnorm[1]{\norm{#1}_{H^s([0,1])}}
\def\int{\intop\limits}
\renewcommand\labelenumi{(\roman{enumi})}
\newcommand\lnorm[1]{\norm{#1}_{\ell^2(\Z)}}
\newcommand\Lnorm[1]{\norm{#1}_{L_2([0,1])}}
\newcommand\LR{{L_2(\R)}}
\newcommand\LRnorm[1]{\norm{#1}_\LR}
\newcommand\Matrix[2]{\hphantom{#1}_#2#1}
\newcommand\norm[1]{\left\|#1\right\|}
\newcommand\ogauss[1]{\left\lceil#1\right\rceil}
\newcommand{\QED}{\hfill
\raisebox{-2pt}{\rule{5.6pt}{8pt}\rule{4pt}{0pt}}%
  \smallskip\par}
\newcommand\Rscalar[1]{\scalar{#1}_\R}
\newcommand\scalar[1]{\left(#1\right)}
\newcommand\Scalar[1]{\scalar{#1}_{[0,1]}}
\newcommand\Span{\mathop{\rm span}}
\newcommand\supp{\mathop{\rm supp}}
\newcommand\ugauss[1]{\left\lfloor#1\right\rfloor}
\newcommand\with{\, : \,}
\newcommand\Null{{\bf 0}}
\newcommand\bA{{\bf A}}
\newcommand\bB{{\bf B}}
\newcommand\bR{{\bf R}}
\newcommand\bD{{\bf D}}
\newcommand\bE{{\bf E}}
\newcommand\bF{{\bf F}}
\newcommand\bH{{\bf H}}
\newcommand\bU{{\bf U}}
\newcommand\cH{{\cal H}}
\newcommand\sinc{{\rm sinc}}
\def\enorm#1{| \! | \! | #1 | \! | \! |}

\newcommand{\dm}{\frac{d-1}{d}}

\let\bm\bf
\newcommand{\bbeta}{{\mbox{\boldmath$\beta$}}}
\newcommand{\bal}{{\mbox{\boldmath$\alpha$}}}
\newcommand{\bbi}{{\bm i}}

\def\nnew{\color{Red}}
\def\mnew{\color{Blue}}

\newcommand{\dI}{\Delta}
\newcommand\aconv{\mathop{\rm absconv}}

\maketitle
\date{}
\begin{abstract}{  
This paper considers the Dirichlet problem
$$
-\div(a\nabla u_a)=f \quad \hbox{on}\,\,\ D, \qquad u_a=0\quad \hbox{on}\,\,\partial D,
$$
for  a Lipschitz domain $D\subset \R^d$, where $a$ is a scalar diffusion function. 
For a fixed $f$, we discuss
under which conditions is  $a$  uniquely determined and  when can $a$ be stably recovered from the knowledge of $u_a$.
 A  first result is that whenever $a\in H^1(D)$, with $0<\lambda \le a\le \Lambda$ on $D$, and   $f\in L_\infty(D)$  is strictly positive,  then 
$$
\|a-b\|_{L_2(D)}\le C\|u_a-u_b\|_{H_0^1(D)}^{1/6}.
$$
More generally, it is shown that   the assumption   $a\in H^1(D)$  can be weakened to      $a\in H^s(D)$, for certain  $s<1$,   at the expense of lowering the exponent
$1/6$  to a   value   that depends on $s$.}

\end{abstract}


\section{Introduction}\label{Introduction}
Let $D$ be a bounded domain (open, connected set) in  $\R^d$, $d \geq 2$.   
We assume throughout the paper that, at a minimum,   $D$ is Lipschitz.   
We define the set of scalar diffusion coefficients
\be
\cA:=\left\lbrace a \in L_\infty(D) \ : \  \lambda\le a\le \Lambda \right\rbrace,
\ee
where $\lambda,\Lambda$ are fixed positive constants.
For $f \in H^{-1}(D)$ (the dual of $H^1_0(D)$) and $a \in \mathcal A$, we consider the elliptic problem
\be
-\div(a\nabla u_a)=f \quad \hbox{on}\,\,\ D, \qquad u_a=0\quad \hbox{on}\,\,\partial D,
\label{ellip}
\ee
written in the usual weak form: $u_a\in H^1_0(D)$ is such that
\be 
\label{i:equation}
\int_D a\nabla u_a \cdot \nabla v =\<f,v\>_{H^{-1(D)},H^1_0(D)}, \qquad  v \in H^1_0(D).
 \ee
Here $H^1_0(D)$ is equipped with the norm $\| v \|_{H^1_0(D)} = \| \nabla v \|_{L_2(D)}$.
The Lax-Milgram theory guarantees that there is a unique solution $u_a \in H^1_0(D)$ of the above problem.

The main interest of the present paper is to understand, for a given $f$, the conditions under  
which the diffusion coefficient $a$ is uniquely determined from the solution $u_a$ to \iref{i:equation}, 
and if so, whether  $a$  can be stably recovered if $u_a$ is known.   
After having fixed $f$, we systematically denote by $u_a$ the solution of  \iref{i:equation}. We
are therefore interested in the stable inversion of the map 
\be
\label{i:smap}
a\mapsto u_a
\ee
which acts from $\cA$ to $H^1_0(D)$.  By stability, we mean 
that when   $u_b$ is close to $u_a$, say in the $H_0^1(D)$ norm, then it follows that $b$ is close to $a$
in some appropriate $L_p(D)$ norm.   The results of this paper will prove such stable inversion but only when certain restrictions are placed on the right side $f$ and further only when the map \eref{i:smap} is restricted to certain subclasses  of $\cA$.

Problems of this type are referred to as  parameter
estimation, or  the identifiability problem in the inverse problems literature, see e.g. \cite{C,A, K, KL,Kn} and the references therein. 
Parameter estimation/identification  for elliptic partial differential equations and their numerical recovery from the (partial) knowledge of $u_a$
is an extensively studied subject that
has been formulated in several settings.  Examples of such settings are   the identifiability of the diffusion 
coefficient $a$ in the problem $-{\rm div}(a\nabla u)=0$ 
from the Neumann boundary data $g$ on $\partial D$, see \cite{KV}, or the  recovery of $a$ from the solution $u$ to 
equation \eref{ellip} supplemented by Dirichlet boundary data, see \cite{Kn}.

Let us make a few elementary remarks about the Dirichlet boundary data 
setting   studied here.  These remarks   extend to other settings as well. 
For $a\in \cA$, we denote by $T_a$ the elliptic operator
$u\mapsto -\div(a\nabla u)$ which is an isomorphism from $H^1_0(D)$ to $H^{-1}(D)$, 
and by $S_a$ its inverse. Then, it is not difficult to check,
see Lemma \ref{lemmabilip} in \S 2, that the map
$a\mapsto S_a$ is bi-Lipschitz from $L_\infty(D)$ to $\cL(H^{-1}(D),H^1_0(D))$, with bounds
\be
\lambda^{2}\|S_a-S_b\|_{\cL(H^{-1}(D),H^1_0(D))} \leq \|a-b\|_{L_\infty(D)}\leq \Lambda^{2}\|S_a-S_b\|_{\cL(H^{-1}(D),H^1_0(D))}, \quad a,b\in\cA.
\label{idS}
\ee
Therefore,  any $a\in \cA$ can be stably identified in the $L_\infty$ norm from the inverse operator $S_a$, that
is, if we knew the solution to \eref{i:equation} for {\it all} possible right  sides then $a$ is uniquely determined. Note that \iref{idS} also means that,
for any $a,b\in\cA$, there exists a right side $f=f(a,b)$, with  $\|f\|_{H^{-1}(D)}=1$, for which we have
the Lipschitz  bound
\be
\label{ab}
\|a-b\|_{L_\infty(D)} \leq  \Lambda^2\|u_a-u_b\|_{H^1_0(D)}.
\ee
 
The $f$ for which  \eref{ab} holds depends on $a$ and $b$. Our objective is  to fix one right  side $f$ 
and study the stable identifiability of $a$ from $u_a$.
It is well known that identifiabiliy cannot hold for an arbitrary right side $f$, even
when $f$ is smooth. For example,  if $u$ is any function in $H^1_0(D)$ such that 
$\nabla u$ is identically $0$ on an open set $D_0\subset D$, then 
setting $f=-\div(a\nabla u)$ for some fixed $a\in \cA$, we find that $u=u_a=u_b$ for any
$b\in\cA$ which agrees with $a$ on $D\setminus D_0$. The above example
can be avoided by assuming that $f$ is strictly positive.  However, even in the case that $f$ is strictly positive, we do not
know a proof of identifiabilty under the general assumption that $a\in\cA$, except in the 
univariate setting.

 In this paper, we show that for strictly positive $f\in L_\infty(D)$,  identifiability and stability hold, for a certain range of $s>0$, in the restricted 
classes $\cA_s\subset \cA$ ,  where
\be
\cA_s:=\cA_{s,M}:=\{a\in\cA\; : \; \|a\|_{H^{s}(D)}\leq M\}.
\label{As}
\ee
Here, $M>0$ is arbitrary but enters in the value of the stability constants. Under such conditions,
we establish  results of the form (see  for example  \ref{cor4.4})
\be
\|a-b\|_{L_2(D)}\leq C\|u_a-u_b\|^{\alpha}_{H^1_0(D)}, \quad a, b\in\cA_s,
\label{genest}
\ee
where the exponent $0<\alpha<1$ depends on $s$
and the constant $C$ depends on $\lambda,\Lambda,\alpha,M,D,f$. Some elementary
observations in the univariate case, see \S 6, show that when $f=1$ and $\cA_s$ includes discontinuous functions, the exponent $\alpha$
cannot be larger than $1/3$.    

There are several existing approaches to establish identifiability.   For the most part, they are  developed 
  for the Neumann problem
\begin{eqnarray}
\label{np}
-{\rm div}(a\nabla u_a)=f \quad \hbox{on}\,\, D, \quad a\frac{\partial u_a}{\partial n}=g \quad \hbox{on}\,\,\partial D,
\end{eqnarray}
where $n$ denotes the outward pointing normal to $\partial D$.
Some approaches use 
singular perturbation arguments, see  \cite{AS}, or the long time behavior of the corresponding unsteady equations, see \cite{HS}.
Some results rely on  the observation that  once $u=u_a$ is given, \eref{np}
may be viewed as a transport equation for the diffusion $a$, see \cite{R1,R2}, and   the  identifiability of 
$a$ from $u_a$ is proven under the assumptions that  $a$ is prescribed 
on the inflow boundary  (the portion of the boundary where $\frac{\partial u_a}{\partial n}<0$) and 
\be
\label{cond}
\inf_D \max\{|\nabla u_a|,\Delta u_a\}>0.
\ee

Other approaches to identifiability  use variational methods, see \cite{KL},  or least-squares techniques, 
see \cite{FP,KS,KW,F}.
These approaches impose strong regularity assumptions on $a$ and $u_a$ as well as the
 assumption 
\be
\label{cond1}
\nabla u_a\cdot \tau>0,
\ee
for a given $\tau \in \mathbb R^d$,
or   the less restrictive condition \eref{cond}.
Rather than directly proving a stability estimate, they derive numerical methods for actually finding the 
diffusion coefficient $a$ from the solution $u_a$ over triangulation $\cT_h$ of $D$ with mesh size $h$. 
One typical reconstruction estimate, see Theorem 1 in \cite{F},  is the following.   Let $r\ge 1$ and let $A_h$ and $V_h$ be the sets of continuous piecewise polynomials 
on $\cT_h$ of degree $r$ and $r+1$, respectively.
If \iref{cond1} holds, and if
$u_a\in  W^{r+3}(L_\infty(D))$ and $a\in H^{r+1}(D)$, then
\be
\label{est}
  \|a-a_h\|_{L_2(D)} \le C \left( h^r + \|u_a - u^{ob}\|_{L_2(D)} h^{-2} \right),
\ee
where $u^{ob}\in L_2(D)$ is an observation of $u_a$, and $a_h\in A_h$ is a
numerical reconstruction of $a$ via least squares type approach from the observation $u^{ob}$.
As shown in Remark \ref{Remark:Falk}, the inequality \eref{est} leads to 
a stability estimate of the form
\be
\label{pip}
\|a-b\|_{L_2(D)}\leq C\|u_a-u_b\|^{\alpha}_{L_2(D)}, \quad \alpha:=\frac{r}{r+2}, \quad a,b\in  \cA_{r+1},
\ee
whenever in addition $u_a,u_b\in  W^{r+3}(L_\infty(D))$ and condition \eref{cond1} holds. Note that $\alpha$ approaches
$1$ as $r\to \infty$.

In summary, the majority of the existing stability estimates are derived for
 solutions to the Neumann problem \eref{np}. As illustrated by \iref{pip},
 they rely on strong regularity assumptions on the 
diffusion coefficients $a$ and on the solutions $u_a$, as well as conditions
on $u_a$ such as  \eref{cond1} or \eref{cond}.   However, one should note that  high order smoothness of $u_a$ 
  generally does not hold, even for smooth $a$ and $f$,  when the domain $D$ does
not have  a smooth boundary.

In this paper, we pursue a variational approach, where we use appropriate  test functions $v$ in \eref{i:equation} 
to derive continuous dependence estimates. We combine these with known elliptic regularity results and obtain 
direct comparison between $\|a-b\|_{L_2(D)}$ and $\|\nabla u_a-\nabla u_b\|_{L_2(D)}$ under milder smoothness
assumptions for the diffusion coefficient $a$, the  domain $D$, and on the right   side $f$, and with {\it no additional 
smoothness assumptions} on $u_a$ and no conditions such as  \eref{cond} or \eref{cond1}.

We mention two special 
cases of our results.  The first, see Corollary \ref{3:cor},  
says that if $D$ is an arbitrary Lipschitz domain, then
for any $f\in L_\infty(D)$ satisfying $f\ge c_f >0$ on $D$, we have   the stability
bound
%
\be
\label{begin}
\|a-b\|_{L_2(D)}\le C \|u_a-u_b\|^{1/6}_{H_0^1(D)}, \quad a,b\in \cA_1.
\ee
We can weaken the smoothness assumption  to the classes $\cA_s$, for $s<1$.    We have two types of results.   In Corollary \ref{vmores}, we prove estimates of the
form 
\be
\label{1:fe}
\|a-b\|_{L_2(D)}\le C \|u_a-u_b\|^{\alpha}_{H_0^1(D)}, \quad a,b\in \cA_s,
\ee
with $\alpha$ depending on $s$, for all  $1/2<s<1$ under the additional assumption that the diffusion coefficients are in VMO and the domain $D$ is $C^1$.  
In Corollary \ref{cor4.4},  we
prove for a general Lipschitz domain $D$,  that  \eref{1:fe} holds for a certain range of $s^*<s<1$ where   we do not require the diffusion coefficients are in VMO but now $s^*$ depends on
properties of the domain $D$.

Estimates  like \eref{pip} have a weaker norm on the right side then those in our results.   However, let us remark that  any   such  estimate 
 can be transformed into an estimate between $\|a-b\|_{L_2(D)}$ and $\|u_a-u_b\|_{L_2(D)}$, if the solutions $u_a$ and $u_b$ 
have more regularity  such as the condition $u_a$ and $u_b$ belong to $H^{1+t}(D)$ for some $t>0$.   For this, one uses the interpolation inequality  
\be
\|v\|_{H^1(D)}\leq C\|v\|_{L_2(D)}^\theta \|v\|_{H^{1+t}(D)}^{1-\theta},\quad v\in H^{1+t}(D), 
\ee
where   $\theta:=\frac t{1+t}$ and
$C_0$ depends only on $D$ and $t$.
 Hence, under the assumption that  $u_a,u_b\in H^{1+t}(D) $,  taking $v=u_a-u_b$, we obtain 
\be
\|u_a-u_b\|_{H_0^1(D)}\leq C_0 \max\{ \|u_a\|_{H^{1+t}(D)}, \|u_b\|_{H^{1+t}(D)}\}^{1-\theta} \|u_a-u_b\|^{\theta}_{L_2(D)},
\ee
which combined with \eref{1:fe} leads to 
\be
\label{weaker}
\|a-b\|_{L_2(D)}\le C \|u_a-u_b\|^{\alpha\theta}_{L_2(D)}.
\ee
Here $C$ depends on the constant in \eref{1:fe},  $C_0$, and  $\max\{ \|u_a\|_{H^{1+t}(D)}, \|u_b\|_{H^{1+t}(D)}\}^{\frac{\alpha}{1+t}}$.

Let us additionally note that  as $r\to\infty$, the result in \eref{pip} leads to  better exponents   then in our results.   This   is caused, at least in part, by the fact that our starting point is \eref{begin} which does not use higher smoothness than $a,b \in H^1(D)$.

Our  paper is organized as follows. In \S\ref{2:section}, we use a variational approach to 
establish a weighted $L_2$ estimate
\be
\|a-b\|_{L_2(w,D)}\leq C\|u_a-u_b\|^{1/2}_{H_0^1(D)},\quad a,b\in \cA_1,
\ee
where the weight is given by $w=a |\nabla u_a|^2 + f u_a$. In order to remove the weight in the above estimate, 
in \S\ref{3:improvements}, we introduce the {\it positivity condition}
\be
 \mbox{\lb{\beta}:} \quad \quad \quad a |\nabla u_a(x)|^2+f(x)u_a(x)\geq c \,{\rm dist}(x,\partial D)^\beta,\quad \hbox{a.e \,\,on}\,\, D, 
\ee
for some $\beta\ge 0$ and $c>0$, 
see Definition~\ref{d:pc}. Under this condition, we prove the stability estimate
\be
\|a-b\|_{L_2(D)}\leq C\|u_a-u_b\|^{\alpha}_{H_0^1(D)},\quad \alpha=\frac{1}{2(\beta+1)}, \quad a,b\in\cA_1.
\ee
Notice that the smaller the $\beta$, the stronger the stability estimate. 

We go further in \S\ref{3:improvements}
 and investigate  which regularity assumptions guarantee that the positivity condition {\bf PC($\beta)$} holds, and thereby obtain results in which this condition is not
 assumed but rather implied by the regularity assumptions on $a$. 
In particular, we
prove that condition \lb{2} is valid for the entire class  $a\in\cA$,  provided $f\in L_2(D)$ with $f\geq c_f>0$.
We also show that    certain smoothness conditions on the diffusion coefficient $a$, the
right side $f$, and  the domain $D$ imply the positivity condition \lb{0}. However, as discussed in \S\ref{ss:necessity_smooth_domain}, 
\lb{\beta} does not generally hold for $\beta<2$ without additional regularity assumptions on the domain $D$. 

In \S\ref{5:section}, we use interpolation arguments to  obtain results under weaker assumptions than  $a,b\in  \cA_1$. 
In \S \ref{s:pwc}, we provide stability estimates in the case when $a$ is piecewise constant
which is not covered by  our general stability results. Finally, in \S 6,  we provide stability estimates in the one dimensional case
for $f=1$ and general $a,b\in \cA$. In this simple case, we also establish converse
estimates which show that the H\"older exponent $\alpha$ in \iref{genest} cannot be above the value $\frac 1 3$
when $a$ and $b$ have low smoothness.

We conclude this introduction by stating some natural open problems
in relation with this paper:
\begin{itemize}
\item[(i)]
While the identifiability problem  is solved in this paper under mild
regularity assumptions, it is still not known whether there exists an $f$ for which the mapping
$a\mapsto u_a$ is injective from $\cA$ to $H^1_0(D)$ for a general multivariate Lipschitz domain $D$.
 \item[(ii)]
The best possible value $\alpha^*=\alpha^*(s)$ of the exponent $\alpha$ in \iref{genest}
is generally unknown. In particular, we do not know if there exists some finite $s_0$
such that $\alpha^*(s)=1$ when $s\geq s_0$.
\item[(iii)]
All our results are confined to the case of scalar diffusion coefficients.
Similar stability estimates for matricial coefficients would require considering
the solutions $u_a$ and $u_b$ for more than one right side $f$.   However 
we are not aware of results that solve this question.
\end{itemize}
 

\section{First estimates}
 \label{2:section}
 
 We begin by briefly discussing the stability properties of the maps
$a\mapsto T_a$ and $a\mapsto S_a$.

\begin{lemma}
\label{lemmabilip}
For any $a,b\in \cA$, we have  
\be
\|T_a-T_b\|_{\cL(H^1_0(D),H^{-1}(D))} =\|a-b\|_{L_\infty(D)},
\label{idT}
\ee
and 
\be
\lambda^{2}\|S_a-S_b\|_{\cL(H^{-1}(D),H^1_0(D))} \leq \|a-b\|_{L_\infty(D)}\leq \Lambda^{2}\|S_a-S_b\|_{\cL(H^{-1}(D),H^1_0(D))}.
\label{idS1}
\ee
\end{lemma}

\noindent
{\bf Proof:} For the proof of \iref{idT}, we observe on the one hand that 
\be
|\<(T_a-T_b)u,v\>_{H^{-1}(D),H^1_0(D)}| \leq \|a-b\|_{L_\infty(D)} \|u\|_{H^1_0(D)}\|v\|_{H^1_0(D)}, \quad u,v\in H^1_0(D),
\ee
which shows that the right quantity dominates the left one in \iref{idT}. On the other hand, 
for any $x\in D$ and $\e>0$ small enough so that the open ball
$B(x,\e)$ of radius $\e$ centered at $x$ is a subset  of $D$, we consider the function $u=u_{x,\e}$ defined by
\be
u(y)=\max\{0,1-\e^{-1}|x-y|\}.
\ee
For such a function, we find that
\be
\<(T_a-T_b)u,u\>_{H^{-1}(D),H^1_0(D)}=C_{x,\e}  \|u\|_{H^1_0(D)}^2,
\quad C_{x,\e}:= |B(x,\e)|^{-1} \int_{B(x,\e)}(a(y)-b(y)) dy.
\ee
By Lebesgue theorem, this shows that 
\be
\|T_a-T_b\|_{\cL(H^1_0(D),H^{-1}(D))}\geq a(x)-b(x), \quad {\rm a.e.} \; x\in D.
\ee
Since we can interchange the role of $a$ and $b$, this 
shows that the left quantity dominates the right one in \iref{idT}. 
For the proof of \iref{idS1}, we observe that $T_a(S_a-S_b)T_b=T_b-T_a$, which yields
\be
\lambda^{2}\|S_a-S_b\|_{\cL(H^{-1}(D),H^1_0(D))}  
\leq \|T_a-T_b\|_{\cL(H^1_0(D),H^{-1}(D))}  \leq \Lambda^{2}\|S_a-S_b\|_{\cL(H^{-1}(D),H^1_0(D))}, \quad a,b\in\cA.
\ee
Combined with \iref{idT}, this gives \iref{idS1}. 
\hfill $\Box$
\newline

As observed in the introduction, the above result does not meet
our objective, since we want to  fix the right   side $f\in H^{-1}(D)$ 
and then study the stable identifiability of $a$ from $u_a$ for all $a\in\cA$.
For such an $f$, let $u_{a},u_{b}$ be the two corresponding 
solutions to \eref{i:equation}, for $a,b\in \cA$. We  use the  notation 
$$
\delta:=a-b,\qquad E:=u_{a}-u_{b}
$$
throughout the paper and we define the linear functional $L: H_0^1(D)\rightarrow \R$, 
$$
L(v):=  \int_D  \delta \nabla u_a\cdot\nabla v,\quad v\in H_0^1(D).
$$
%
By subtracting the two weak equations \eref{i:equation} for $a$ and $b$,  we derive another representation of $L$, 
\be
\label{2:diff1}
L(v)= -\int_D b\nabla E\cdot \nabla v ,\quad v\in H_0^1(D).
\ee	
The following theorem gives two basic estimates for bounding the difference $\delta=a-b$.
The first one illustrates that difficulties arise when $a-b$ changes sign, while the second puts forward the role of the weight $w=a|\nabla u_a|^2 + f u_a$.

\begin{theorem} 
\label{2:theorem1}  Let $D$ be a Lipschitz domain.  Consider equation \eref{i:equation} with diffusion coefficients $a$ and $b$.
The  following two inequalities hold for $\delta:=a-b$.  
\vskip .1in

\noindent
{\rm (i)} For  any $a,b\in\cA$ and $f\in H^{-1}(D)$, we have
$$
\left|\int_D \delta |\nabla u_a|^2 \right| \le \Lambda \|f\|_{H^{-1}(D)}   \|E\|_{H_0^1(D)}.
$$

\noindent
 {\rm (ii)} For any  $a,b\in \cA_1$ and   $f\in L_\infty(D)$, we have  
\be
 \label{2:estimate2}
  \int_D\frac{\delta^2}{a^2}\(a |\nabla u_a|^2   + f u_a\) \le C_0\|E\|_{H_0^1(D)}, 
\ee
where
\be
\label{2:defC} 
 C_0:=C\|f\|_{L_\infty(D)}(1+\max\{\|\nabla a \|_{L_2(D)}, \|\nabla b \|_{L_2(D)}\}),
 \ee
and $C$ is a constant depending only on $D,d,\lambda,\Lambda$.
 \end{theorem}

\noindent
{\bf Proof:}  To prove  (i), we take  $v=u_a\in H_0^1(D)$ and  obtain
$$
L(u_a)=  \int_D \delta |\nabla u_a|^2.
$$
Using this in \eref{2:diff1} yields
\be
\label{2:estimate12}
 \int_D \delta |\nabla u_a|^2= -   \int_D b\nabla E \cdot \nabla u_a \le \Lambda \|u_a\|_{H_0^1(D)}\|E\|_{H_0^1(D)}.
 \ee
If we take $v=-u_a$, we derive the same estimate for the negative of the left side of \eref{2:estimate12} which yields (i).

To prove (ii), we  define $\bar \delta:=\delta/a$ which belongs to $H^1(D)$ since $a,b\in\cA_1$.  
Integrating by parts, we 
have  for any $v\in H_0^1(D)$,  
 \be
 \label{2:deltabar} 
L(v)= \int_D\bar \delta a \nabla u_a \cdot \nabla v  =  - 
\int_D \nabla \bar\delta \cdot \nabla u_a a v - \int_D \bar \delta\div (a \nabla u_a)v.
\ee
Since  $f=-\div (a \nabla u_a)$, this gives
\be
\label{2:deltabar1}
L(v)= \frac 1 2 \int_D  \bar \delta a \nabla u_a \cdot \nabla v - \frac 1 2  \int_D \nabla  \bar\delta \cdot \nabla u_a av 
 + \frac 1 2 \int_D \bar \delta f  v, \quad  v\in H_0^1(D).
\ee
Now, we chose  $v = \bar\delta u_a\in H_0^1(D)$  to obtain
\be
\label{2:deltabar2}
L(\bar \delta u_a)= \frac 1 2 \int_D \bar\delta^2 a |\nabla u_a|^2  
 + \frac 1 2 \int_D \bar \delta^2 f u_a.
\ee
Inserting \eref{2:deltabar2} into \eref{2:diff1} results in
\be
\label{2:deltabar3}
 \frac 1 2 \int_D \bar\delta^2 a |\nabla u_a|^2  
 + \frac 1 2 \int_D \bar \delta^2 f u_a 
= 
 -\int_D b \nabla E \cdot \nabla (\bar \delta u_a)\le \Lambda \|\nabla(\bar \delta u_a)\|_{L_2(D)}\|E\|_{H_0^1(D)}.
\ee
Now, we resort to the estimate (see e.g.  Chapter 8 in  \cite {GT})
$$
\|u_a\|_{L_\infty(D)}\leq C\|f\|_{L_\infty(D)},
$$ 
 where $C$ depends only on $\lambda,\Lambda$ and $D$ (throughout the rest of this proof $C>0$ will be  a generic constant that depends on at most  $d,D,\lambda,\Lambda$).
 We use this result  together with the energy estimate
 $$
 \| \nabla u_a \|_{L_2(D)} \leq \| f \|_{H^{-1}(D)} \leq C \| f \|_{L_\infty(D)}
 $$
 to  obtain  the bound
\begin{eqnarray} 
\label{2:deltabar4}
\|\nabla(\bar \delta u_a)\|_{L_2(D)}&\le& \left\|\frac{\delta}{a}\right\|_{L_\infty(D)}\|\nabla u_a\|_{L_2(D)}+\left\|\frac{u_a}{a}\right\|_{L_\infty(D)}\|\nabla \delta\|_{L_2(D)}+
\left\|\frac{\delta}{a^2}\right\|_{L_\infty(D)}\left\|u_a\right\|_{L_\infty(D)}\|\nabla a\|_{L_2(D)}\nonumber\\
&\le & 2\Lambda \lambda^{-1}\|\nabla u_a\|_{L_2(D)}+ 
\lambda^{-1}\|u_a\|_{L_\infty(D)}\| \nabla \delta\|_{L_2(D)} +      2 \Lambda \lambda^{-2}\|u_a\|_{L_\infty(D)}\|\nabla a\|_{L_2(D)}\nonumber \\
&\le & C\|f\|_{L_\infty(D)}(1+\max\{\|\nabla a \|_{L_2(D)}, \|\nabla b \|_{L_2(D)}\}).
\end{eqnarray}
Finally, plugging this estimate into \eref{2:deltabar3}, we derive that 
 \begin{eqnarray} 
 \nonumber
  \int_D \frac{\delta^2}{ a} |\nabla u_a|^2  
 +  \int_D \frac{\delta^2}{a^2} f u_a
&=& \int_D \bar\delta^2 a |\nabla u_a|^2  
 +  \int_D \bar \delta^2 f u_a 
\le 2\Lambda \|\nabla(\bar \delta u_a)\|_{L_2(D)}\|E\|_{H_0^1(D)}\\
\nonumber
&\leq&C\|f\|_{L_\infty(D)}(1+\max\{\|\nabla a \|_{L_2(D)}, \|\nabla b \|_{L_2(D)}\} )\|E\|_{H_0^1(D)},
\end{eqnarray}
and the proof is completed.
 \hfill $\Box$
 \newline
 
 Note that  when $a\leq b$ or $b\leq a$ a.e. on $D$ and condition \eref{cond1} holds in the sense that
$\nabla u_a\cdot \tau\geq c>0$, then part {\rm (i)} gives the stability estimate
 $$
 \|a-b\|_{L_1(D)}\leq C\|f\|_{H^{-1}(D)}\|u_a-u_b\|_{H_0^1(D)}.
 $$
 However, we can not claim such a result if the difference $(a-b)$ changes sign on a subset of $D$ with a positive measure. 
 In the sequel of the paper, we will not use (i), and instead rely only on (ii).
  
 
\section{Improvements of Theorem \ref{2:theorem1}}
\label{3:improvements}

Theorem \ref{2:theorem1} is not satisfactory as it stands, since we want to replace the left side of \eref{2:estimate2}, 
by $\|a-b\|_{L_2(D)}^2$. Obviously, this is possible when there exists a constant $c >0$ such that the weight
satisfies
\be
\label{pa}
a|\nabla u_a |^2+fu_a\ge c \quad \hbox{a.e. \,\,on}\quad D.
\ee
In order to understand this condition, suppose that $f$ does not change sign. In that case, the weak maximum principle \cite{GT} guarantees that $u_a$ has the same sign as $f$ and therefore the product $u_a f \geq 0$. 
Hence, \eqref{pa} requires that $u_a$ and $|\nabla u_a|$ do not vanish simultaneously. 
We   prove in \S\ref{beta0}  that such a constant $c$ exists provided certain   (strong)   smoothness assumptions  for the diffusion coefficient $a$, the right side $f$, and the domain $D$ hold.  
However, in order to allow milder regularity assumptions, we introduce the following weaker positivity condition.

\begin{definition}[Positivity Condition]\label{d:pc}
We say that $(D,f,a)$ satisfy the positivity condition \lb{\beta} if there exists a constant $c>0$ such that 
\be
\label{3:a2}
a(x)|\nabla u_a(x)|^2+f(x)u_a(x)\ge c \dist(x,\partial D)^\beta, \quad \hbox{a.e.}\  x\in D.
\ee
\end{definition}

Notice the positivity condition \lb{0} is \eref{pa}. 
In  Lemma \ref{Gr}, we show 
that for every Lipschitz domain $D$ and $a\in \cA$, we have that $(D,a,f)$ satisfies the positivity condition \lb{2} provided $f$ is strictly positive and in $L_2(D)$.
In fact, in this case, the  constant $c$ in \iref{3:a2} is uniform over the class $\cA$.
In addition, we provide examples which show that additional regularity 
assumptions are required for $(D,a,f)$ to satisfy the positivity condition \lb{\beta}
if $\beta<2$.
For now, we prove the following theorem which shows how a positivity condition \lb{\beta} guarantees a stability estimate of the type we want.  

\begin{theorem} 
\label{3:inversetheorem}
 Let $D$ be a Lipschitz domain.
 Assume that $a$, $b\in \cA_1$, $f\in L_\infty(D)$ and denote by $u_a$, $u_b$ 
 the corresponding solutions to {\rm\eref{i:equation}}. If $(D,a,f)$ satisfies the positivity condition \lb{\beta} for $\beta\geq 0$, then we have

\be
\label{3:it}
\|a-b\|_{L_2(D)}\le C\sqrt{1+C_0} \|u_a-u_b\|_{H_0^1(D)}^{\frac{1}{2(\beta+1)}} ,
 \ee  
 where $C_0$ is the constant from {\rm \eref{2:defC}} and $C$ is
 a constant depending only on $D,d,\lambda,\Lambda$, and $c$ the constant in \eqref{3:a2}.
\end{theorem}
\noindent
{\bf Proof:}
We recall the notation $\delta=a-b$, $E=u_a-u_b$, and  start with  the weighted $L_2$ estimate \eqref{2:estimate2} provided in Theorem~\ref{2:theorem1}, namely
\begin{equation}\label{e:starting_point}
\int_D \frac{\delta^2}{a^2} w \leq C_0 \| E \|_{H^1_0(D)},  \quad w:=  a  | \nabla u_a |^2+ fu_a,
\end{equation}
where $C_0$ is the constant in \eref{2:defC}.   This proves the result in the case $\| E \|_{H^1_0(D)}=0$ since $w>0$ on $D$.   Therefore, in going further, we assume
$\| E \|_{H^1_0(D)}>0$.

 The presence of the non-negative weight $w$ is handled by decomposing the domain $D$ into two sets 
$$
D_\rho := \{ x \in D \ : \ \dist(x,\partial D) \geq \rho \} \qquad \textrm{and} \qquad D^c_\rho := D \setminus D_\rho,
$$
where $\rho>0$ is to be chosen later. 
The triplet $(D,a,f)$ satisfies the positivity condition \lb{\beta}, which guarantees that $w  \geq c \rho^\beta$ on $D_\rho$.
Hence, we deduce  that
\begin{equation}\label{e:Drho}
\int_{D_{\rho}} \delta^2  \leq  \Lambda^2c^{-1} \rho^{-\beta}\int_D \frac{ \delta^2}{a^2}w\le 
\Lambda^2c^{-1}C_0  \rho^{-\beta} \| E \|_{H^1_0(D)}.
\end{equation}
On $D^c_\rho$, the Lipschitz regularity assumption on $\partial D$ implies the existence of a constant $B$ such that $|D^c_\rho|\le B\rho$. 
As a consequence, we obtain
\begin{equation}\label{e:DrhoC}
\int_{D^c_{\rho}} \delta^2  \leq 4\Lambda^2 | D^c_\rho| \leq  4\Lambda^2 B\rho.
\end{equation}
Combining the last two estimates with the choice $\rho = \| E \|_{H^1_0(D)}^{\frac{1}{\beta+1}}$ proves \eref{3:it} and ends the proof. \hfill $\Box$

\subsection {The positivity condition \lb{0}}
\label{beta0}
In view of the exponent in \eref{3:it}, the strongest stability occurs when $\beta=0$.  In this section, we show that if $(D,a,f)$ are sufficiently smooth then  \lb{0} is satisfied.
We denote by $C^{k,\alpha}(D)$, $k\in \N_0$, $0<\alpha\leq 1$,  the H\"older spaces equipped with the semi-norms
$$
|f|_{C^{k,\alpha}(D)}:=\sup_{|\gamma|=k}\sup_{x,y\in D,\,\,x\neq y}\left\{\frac{|\partial^\gamma f(x)-\partial^\gamma f(y)|}{|x-y|^\alpha}\right\},
$$
and norms
$$
\|f\|_{C^{k,\alpha}(D)}:=\sup_{|\gamma|\leq k}\|\partial^\gamma f\|_{L_\infty(D)}+|f|_{C^{k,\alpha}(D)}.
$$

\subsubsection{Sufficient conditions}
The following lemma gives a sufficient condition for $(D,a,f)$ to satisfy the positivity condition \lb{0}.

\begin{lemma}
\label{3:Schtlemma}
 Assume that for some $\alpha>0$,   $D$ is  a $C^{2, \alpha}$  domain and  $f\in C^{0,\alpha}(D)$   with $f\geq c_f>0$.
Furthermore, assume that the diffusion coefficient $a$ belongs to $\cA \cap C^{1,\alpha}(D)$, with
\be
\|a\|_{C^{1,\alpha}(D)}\leq A.
\ee
Then, the triplet $(D,a,f)$ satisfies the positivity condition \lb{0},
with constant $c$ depending on $D,\lambda,\Lambda, \|f\|_{ C^{0,\alpha}}, c_f$ and $A$.
\end{lemma}

\noindent
{\bf Proof:}  
 We have that 
$$
a(x)|\nabla u_a(x)|^2+f(x)u_a(x)\geq \min\{\lambda,c_f\}\left (|\nabla u_a(x)|^2+u_a(x)\right ),
$$
since $u_a\geq 0$ according to the weak maximum principle \cite{GT}.
We proceed by showing that $|\nabla u_a|^2+u_a\geq c$, a.e. on $D$. We do this  by 
contradiction.   Assume that there exists a sequence 
$\{a_n\}_{n\geq 0}$ of diffusion coefficients $a_n\in \cA$ with $\|a_n\|_{C^{1,\alpha}(D)}\leq A$  such that, for each $n\geq 0$, there exists $x_n\in D$ with
\be
|\nabla u_{a_n}(x_n)|^2+u_{a_n}(x_n) \leq \frac 1 n.
\ee
Note that  the  assumptions of the theorem imply that the equation \eref{i:equation} holds in the strong sense. 
Then, the classical Schauder estimates, see \cite{GT}, tell us that
\be
\|u_{a_n}\|_{C^{2,\alpha}(D)}\leq C,
\ee
where $C$ depends on $A$, $D$, $\alpha$, $\lambda$ and $\Lambda$. 
Then by compactness, up to a triple subsequence extraction, we may assume that
\begin{enumerate}
\item
$a_n$ converge in $C^1$ towards a limit $a^*$,
\item
$u_{a_n}$ converges in $C^2$ towards a limit $u^*$,
\item
$x_n$ converges in $\o D$ towards a limit $x^*$.
\end{enumerate}
Therefore, the equation 
\be
-a^*\nabla u^*-\nabla a^* \cdot \nabla u^*=f,
\ee
is satisfied on $D$, with homogeneous boundary conditions, and we have
\be
u^*(x^*)=0\quad {\rm and}\quad \nabla u^*(x^*)=0.
\ee
The first equality shows that $x^*$ lies on the boundary, due to the strong maximum
principle, and therefore the second equality contradicts the Hopf lemma, see \cite{GT}.  
\hfill $\Box$
\newline

We  have the following corollary.

\begin{corollary}
\label{3:Schtheorem}
  Assume that  for some $\alpha>0$, $D$ is a $C^{2,\alpha}$ domain, $f\in C^{0,\alpha}(D)$   with $f\geq c_f>0$ and 
the diffusion coefficient $a \in \cA \cap C^{1,\alpha}(D)$, with $\|a\|_{C^{1,\alpha}(D)}\leq A$ .
Furthermore, assume that $b\in \cA_1$.
Let $u_a$ and $u_b$  be  the corresponding solutions to {\rm \eref{i:equation}},
then
\be
\label{3:Sch}
\|a-b\|_{L_2(D)}\le C_0\|u_a-u_b\|^{1/2}_{H_0^1(D)},
\ee
where $C_0=C\|f\|^{1/2}_{L_\infty(D)}(1+\max\{\|\nabla a \|_{L_2(D)}, \|\nabla b \|_{L_2(D)}\})^{1/2}$
and $C$ is a constant depending only on $D,d,\lambda,\Lambda, c_f,\|f\|_{C^{0,\alpha}}$, and $A$. 
In particular, under the same assumptions on $D$, $f$, and $b$, we have the estimate
\be
\label{3:Schsob}
\|a-b\|_{L_2(D)}\le C_s\|u_a-u_b\|^{1/2}_{H_0^1(D)}, \quad a\in\cA_s,
\ee
for all  $s>1+\frac d 2$.
\end{corollary}

\noindent
{\bf Proof:}   The inequality  \iref{3:Sch}  follows from  Theorem \ref{3:inversetheorem} and Lemma \ref{3:Schtlemma},
while  \iref{3:Schsob} follows by the Sobolev embedding of $H^s$ into the relevant H\"older spaces.
\hfill $\Box$

\subsubsection{ The condition {\bf PC}($\beta$), $\beta<2$, requires smooth domains } 
\label{ss:necessity_smooth_domain}

In this section, we show that we cannot expect the triplet $(D,a,f)$  to satisfy a positivity condition \lb{\beta}, $\beta<2$, without additional regularity assumptions  on the domain $D$.
We consider the problem, 
\begin{eqnarray}
 \label{poisson}
-\Delta u&=&1, \quad \hbox{on}\,\, D=(0,1)^d,
\\ \nonumber
u&=&0,  \quad \hbox{on}\,\,\partial D,
 \end{eqnarray}
corresponding to the case $a=1$, $f=1$, $D=(0,1)^d$.  We begin with    the following lemma.

\begin{lemma}
  \label{poisson-smooth}
  The solution $u$  to \eref{poisson} is  in the H\"{o}lder space $C^{1,\alpha}(D)$
  for all  $0<\alpha <1$.
  \end{lemma}
 
 \noindent
{\bf Proof: }
The solution $u$ can be expanded in the eigenfunction basis
\be \label{e:sol_cube}
u(x)=\sum_{n\in \N^d} c_n s_n(x), \quad s_n(x):= \prod_{i=1}^d\sin\(\pi n_i x_i\),\quad x=(x_1,\ldots,x_d),
\ee
with  coefficients $c_n$, $n=(n_1,\dots,n_d)$, given by the formula
$$
 c_n= \left\{ \begin{array}{ll} \frac {4^{d}}{\pi^{2+d}(n_1^2+\cdots+n_d^2)n_1\dots n_d},& \hbox{if \,\,all}\,\, n_i \,\,\hbox {are \,\, odd}, \\ \\
 0, & \hbox{otherwise}.\end{array} \right.
 $$
To  prove the stated smoothness for the 
partial derivative $\frac{\partial u}{\partial x_1}$, we first show that
\be
\label{derbd}
\sum_{n\in\N^d}\frac {1}{(n_1^2+\cdots+n_d^2)n_2\dots n_d}<\infty.
\ee
For this, we use the fact that, for any $A>0$,
$$
\sum_{k\geq 1} (A+k^2)^{-1}\leq \int_{0}^\infty (A+t^2)^{-1}dt = \frac \pi {2\sqrt A},
$$
and thus
$$
\begin{disarray}{ll}
\sum_{n\in\N^d}\frac {1}{(n_1^2+\cdots+n_d^2)n_2\dots n_d}&  \leq \frac{\pi}{2}
\sum_{(n_2,\dots,n_d)\in\N^{d-1}} \frac 1 {n_2\dots n_d\sqrt{n_2^2+\cdots+n_d^2}} \\\\
& \leq \frac{\pi}{2(d-1)^{\frac12}}
\sum_{(n_2,\dots,n_d)\in\N^{d-1}} \frac 1 {(n_2\dots n_d)^{1+\frac 1 {d-1}} }\\\\
&=\frac \pi {2(d-1)^{\frac12}}\(\sum_{k\geq 1} k^{-1-\frac 1 {d-1}} \)^{d-1}<\infty,
\end{disarray}
$$
where we have used the inequality between the arithmetic and geometric mean of $n_2^2,\ldots,n_d^2$.

From \eref{derbd}, we can differentiate $u$ termwise and obtain that  $\frac{\partial u}{\partial x_1}$ is continuous.
The same holds for all other partial derivatives, and thus $u\in C^1(D)$. In order to prove that 
$u$ belongs to the H\"older space $C^{1,\alpha}(D)$ for sufficiently small $\alpha>0$, it suffices 
to check in addition that
$$
\sum_{n\in\N^d}\frac {n_i^\alpha}{(n_1^2+\cdots+n_d^2)n_2\dots n_d}<\infty, \quad i=1,\dots,d.
$$
Each term in this series is less than $\frac 1{(n_1^2+\cdots+n_d^2)^{1-\frac \alpha 2}n_2\dots n_d}$.
We thus proceed to a similar computation using the fact that
$$
\sum_{k\geq 1} (A+k^2)^{-1+\frac \alpha 2}\leq \frac C{(\sqrt A)^{1-\alpha}},
$$
and derive that
$$
\sum_{n\in\N^d}\frac {n_i^\alpha}{(n_1^2+\cdots+n_d^2)n_2\dots n_d}
\leq C\(\sum_{k\geq 1} k^{-1-\frac {1-\alpha} {d-1}} \)^{d-1}<\infty,
$$
since $\alpha<1$.
\hfill $\Box$
\newline 

The above lemma allows us to show that the positivity condition
\lb{\beta} does not hold for $\beta<2$, and in particular when $\beta=0$ when $D=(0,1)^d$.

\begin{proposition}
  \label{poisson-beta}
  Let $D=(0,1)^d$ and $a=f=1$, with $d\geq 2$. Then the triplet $(D,a,f)$ does not satisfy the positivity condition \lb{\beta}
  if $\beta<2$.
  \end{proposition}
  
  \noindent
{\bf Proof: }
As shown in Lemma \ref{poisson-smooth}, the solution $u$ to \eref{poisson} is in the class
$C^{1,\alpha}(D)$ for all $0<\alpha<1$, and therefore $\nabla u$ can be continuously extended up to the boundary $\partial D$.
 Since the tangential derivatives of $u$ vanish on the boundary, it follows
that when $x^*$ is a corner of the cube $[0,1]^d$, then  $\nabla u(x^*)=0$. By H\"older regularity,
we find that
\be
|\nabla u(x)| \leq  C{\rm dist}(x,x^*)^\alpha \quad {\rm and}\quad |u(x)|\leq C{\rm dist}(x,x^*)^{1+\alpha}, \quad x\in D,
\ee
and therefore
\be
a(x)|\nabla u_a(x)|^2 +f(x)u_a(x)\leq C{\rm dist}(x,x^*)^{2\alpha}, \quad x\in D,
\ee
for all $0<\alpha<1$.
Thus, \lb{\beta} cannot hold for  any $\beta<2$.\hfill $\Box$


\subsection {The positivity condition \lb{2}}
\label{beta2}

In this section, we  show that the triplet $(D,a,f)$ satisfies the positivity condition  \lb{2} for any 
 Lipschitz domain $D$, any $a\in\cA$, and any $f\in L_2(D)$, with $f\geq c_f>0$.    
 For this, we use the lower bounds on the Green functions
established in \cite{GW}. 

\begin{lemma}
\label{Gr}
Let $D$ be a Lipschitz domain,  $a\in\cA$, and $f \in L_2(D)$ with $f\geq c_f>0$. Then the triplet $(D,a,f)$ satisfies the positivity condition \lb{2} with a constant $c$ only depending on $\lambda, \Lambda, d, D, c_f$.
\end{lemma}

\noindent
{\bf Proof:}
In this proof, $C$ denotes a generic constant only depending on $D,\lambda,\Lambda,d,c_f$.
We recall that for every $y \in D$, there exists a unique Green's function $G_a(\cdot, y) \in W^1_{0}(L_1(D))$, such that
$$
  \int_D \nabla G_a(x,y) \nabla v(x) \, dx = v(y), \quad v \in C_0^\infty(D).
$$
One can show that  
$$
G_a(x,y)\geq C|x-y|^{-(d-2)}, \quad \hbox {for}\quad |x-y| \leq \frac 1 2 \rho(x),\quad d\geq 2,
$$
where $\rho(x):=\dist(x,\partial D)$. A proof of this fact in the case $d\geq 3$ can be found in \cite[Theorem 1.1]{GW}. The same proof 
holds also in the case $d=2$,  utilizing  the regularity properties of the two dimensional  Green's function discussed in \cite{DK}.

Now, given any $x\in D$, let $B(x, \rho(x)/2) \subset D$  be the ball centered at $x$ with radius $\rho(x)/2$.
Since $G_a(x,y)\geq 0$, $x,y\in D$,
we have
\begin{eqnarray}
\nonumber
u_a(x)&=&\int_D f(y)G_a(x,y)\, dy \geq \int_B f(y)G_a(x,y)\, dy 
\\ \nonumber
&\ge& C  \int_{B(x,\rho(x)/2)} |x-y|^{-(d-2)}\, dy \ge C \rho^2(x)= C[{\rm dist}(x,\partial D)]^2,
\nonumber
\end{eqnarray}
and the desired result follows.
\hfill $\Box$

We have the following corollary.
\begin{corollary}
\label{3:cor}
 Let  $D$ be a Lipschitz domain, $a,b\in \cA_1$,  $f\in L_\infty(D)$ with 
 $f\geq c_f>0$, and $u_a, u_b \in H^1_0(D)$ 
 be the corresponding solutions to {\rm \eref{i:equation}},
  then we have
\be
\label{3:it1}
\|a-b\|_{L_2(D)}\le  C\sqrt{1+C_0}\|u_a-u_b\|_{H_0^1(D)}^{1/6},
 \ee 
 where {$C_0$ is the constant in {\rm\eref{2:defC}} and} $C$ is a 
 constant depending only on $D,d,\lambda,\Lambda$ and the minimum $c_f$ of $f$.
\end{corollary}
\noindent
{\bf Proof:} 
The proof follows from Theorem \ref{3:inversetheorem} and Lemma \ref{Gr}.
\hfill $\Box$

 \section{Finer estimates for parameter recovery}
 \label{5:section}
 
 We have proved  Corollary \ref{3:cor}  
 for Lipschitz domains $D$ under the assumptions that 
 $a,b\in \cA_1$
 and $f\in L_\infty(D)$, with  $f\geq c_f>0$.   In this section, we shall weaken 
 the smoothness assumption on $a$ and $b$ at the expense of decreasing the exponent 
 $1/6$ appearing on the right   side of \eqref{3:it1}.

 \subsection{Finer estimates}

 Our method for reducing the smoothness assumptions on the diffusion coefficients in the 
 stability Theorem \ref{3:inversetheorem} will be 
based on interpolation.   We recall that if $a\in H^s(D)$, where $D\subset \R^d$ is  a bounded Lipschitz domain,
then   for each $t>0$, there is a function $a_t\in H^1(D)$ satisfying  the inequality 
\be
\label{3:Kf}
\|a-a_t\|_{L_2(D)}+ t\|\nabla a_t\|_{L_2(D)}\le Ct^s\|a\|_{H^s(D)},
\ee
where the constant $C$ depends only on $D$.  
Note that the standard construction of $a_t$ is a local mollification of $a$, and therefore $a_t \in \mathcal A$ whenever $a\in \mathcal A$.

Our stability estimate relies on  the following result which can be derived from Theorem 2.1 in \cite{BDN}:
\begin{lemma} 
Given $a, b \in \mathcal A$, assume that for some $0<\theta \leq 1$ there exists a constant $M$ such that
$$ 
\| \nabla u_a \|_{L_{2/(1-\theta)}(D)} \leq M.
$$ 
Then, 
\begin{equation}
\label{e:theta_estimate}
\| u_a - u_b \|_{H^1_0(D)} \leq \lambda^{-1}(2\Lambda)^{1-\theta}M\|a-b\|^{\theta}_{L_2(D)}.
\end{equation} 
\end{lemma}
\noindent
{\bf Proof:}   We take $p=\frac{2}{1-\theta}$ in Theorem 2.1 of \cite{BDN}, then for $q=\frac{2}{\theta}$,  we have from (2.2) of \cite{BDN}
\be
\label{bdn1}
\| u_a - u_b \|_{H^1_0(D)} \leq \lambda^{-1}M\|a-b\|^{\theta}_{L_q(D)}.
\ee
Since $\|a-b\|_{L_q(D)}\le \|a-b\|_{L_2(D)}^{2/q}(2\Lambda)^{1-2/q}$, the lemma follows.\hfill $\Box$
\vskip .1in

This motivates the following definition.
\begin{definition}[Gradient Condition]
We say that a function $u \in H^1_0(D)$ satisfies the gradient condition \grad{}, $0<\theta\leq 1$, if  
\begin{equation}
\label{e:lpbound}
\| \nabla u \|_{L_{2/(1-\theta)}(D)}  \leq M.
\end{equation}
\end{definition} 

%
%
%
%

 We now prove our main result regarding stable recovery of parameters provided that $u_a$ satisfies the gradient condition \grad{}.
Later, in \S\ref{solution-smoothness},   we elaborate on what classical 
smoothness conditions on the diffusion coefficient $a \in \mathcal A$ guarantees that this gradient condition holds.
%
%

   \begin{theorem}
 \label{5:istheorem1}
 Let $D$ be a Lipschitz domain, $f\in L_\infty(D)$ with $f\geq c_f>0$, and $a,b\in \cA_s$ for some $1/2<s\leq 1$.
Let $u_a, u_b\in H^1_0(D)$ be the corresponding solutions to \eqref{i:equation}.
If  $u_a, u_b$ both satisfy the gradient condition \grad{} for some $\frac{1-s}{s}<\theta\leq 1$,
then we have
  \be
 \label{3:t2}
  \|a-b\|_{L_2(D)}\le C
\sqrt{1+(\|a\|_{H^s(D)}+\|b\|_{H^s(D)})^{\frac {1} {3s}}} \|u_a-u_b\|_{H_0^1(D)}^{\frac{1}{6}-\frac{1-s}{6s\theta}},
 \ee
 where $ C$  is a constant depending only on $D,d,\theta, \lambda,\Lambda$,  the minimum $c_f$ of $f$, $\|f\|_{L_\infty(D)}$, and $M$.
 \end{theorem}
 \noindent
 {\bf Proof:}  
 We use the notation
 $$
 E:=u_a-u_b, \quad E_t:=u_{a_t}-u_{b_t}, \quad \delta:=a-b, \quad \delta_t:=a_t-b_t,
 $$
 where $a_t,b_t\in\cA_1$ are the functions satisfying \eref{3:Kf}.
Throughout  the proof $C>0$ will be a generic constant that depends on at most  
 $D,d,\theta,\lambda,\Lambda$, $M$, $\|f\|_{L_\infty(D)}$, and  the minimum $c_f$ of $f$.
 In what follows, the  value of $C$
may change at each appearance. We denote by 
\begin{equation}\label{e:bound_M} 
M_0:=\|a\|_{H^s(D)}+\|b\|_{H^s(D)}\ge \|a\|_{L_2(D)}+ \|b\|_{L_2(D)}\ge 2\lambda | D|^{1/2}.
\end{equation}
It follows from  \eref{3:Kf} that 
\be
\label{know111}
\|\delta-\delta_t\|_{L_2(D)}\le C M_0t^s.
\ee
We want to  bound $\|\delta\|_{L_2(D)}$.  For this, we define
the set $D_\rho := \{ x \in D \ : \ \dist(x,\partial D) \geq \rho \}$, with the   value of $\rho>0$ to be chosen shortly.
Using \eref{know111}, we find that
\begin{eqnarray}
\label{popi}
\|\delta\|^2_{L_2(D)}&=&\|\delta\|^2_{L_2(D_\rho^c)}+\|\delta\|^2_{L_2(D_\rho)}
\leq \|\delta\|^2_{L_2(D_\rho^c)}+2\|\delta-\delta_t\|^2_{L_2(D)} +2\|\delta_t\|^2_{L_2(D_\rho)}
\nonumber
\\
&\le& \|\delta\|^2_{L_2(D_\rho^c)}+CM_0^2t^{2s}+2\|\delta_t\|^2_{L_2(D_\rho)}.
\end{eqnarray}
To estimate the two norms above, we proceed as in the  proof of Theorem \ref{3:inversetheorem}.
First, for $a$, $b\in \cA$ and a Lipschitz domain $D$ we have
\be
\label{complement}
\|\delta\|^2_{L_2(D_\rho^c)}=\int_{D^c_{\rho}} \delta^2  \leq 4\Lambda^2 | D^c_\rho| \leq C  \rho;
\ee
see \eqref{e:DrhoC}.
Since  $a_t$ and $b_t$  are in $ \cA_1$, according to Lemma~\ref{Gr}, $(D,a_t,f)$ and $(D,b_t, f)$  satisfiy the positivity condition \lb{2} with a constant $c$ only depending on $\lambda,\Lambda, D,d$.  
Hence \eref{e:Drho} holds with $\beta=2$ and therefore, we have
$$
\|\delta_t\|^2_{L_2(D_\rho)}=\int_{D_{\rho}} \delta_t^2  \le C\rho^{-2}(1+ \max\{\|\nabla a_t\|_{L_2(D)}, \|\nabla b_t\|_{L_2(D)}\})\|E_t\|_{H_0^1(D)}.
$$
This, together with \eqref{3:Kf} implies that
\begin{equation}
\|\delta_t\|^2_{L_2(D_\rho)} \leq C\rho^{-2} (1+M_0t^{s-1})\|E_t\|_{H_0^1(D)}.
\label{klm}
\end{equation}
We  substitute \eref{complement} and \eref{klm} into \eref{popi} 
to arrive at
\be
\label{popi1}
\|\delta\|^2_{L_2(D)}\leq C\rho+CM_0^2t^{2s}+C\rho^{-2} (1+M_0t^{s-1})\|E_t\|_{H_0^1(D)}.
\ee

We  now proceed to estimate   $\|E_t\|_{H_0^1(D)}$ by taking advantage of  the gradient condition \grad{} satisfied by $u_a$ and $u_b$.
Since $u_a$ satisfies the gradient condition \grad{} and $a_t\in\cA$, it follows from the stability estimate \eqref{e:theta_estimate} that
\be
\label{know11}
\|u_a-u_{a_t}\|_{H_0^1(D)} \le  C\|a-a_t\|_{L_2(D)}^\theta \le C  (M_0t^s)^\theta.
\ee
The same estimate holds with $a$ replaced by $b$,
and therefore
\be
\label{knowE}
\|E_t\|_{H_0^1(D)}\leq \|u_{a_t}-u_{a}\|_{H_0^1(D)} +\|u_{a}-u_{b}\|_{H_0^1(D)} +\|u_{b}-u_{b_t}\|_{H_0^1(D)} 
\leq C(M_0t^s)^\theta+\|E\|_{H_0^1(D)}.
\ee
Placing this estimate into \eref{popi1}  gives
\be
\label{finally1}
\|\delta\|^2_{L_2(D)}\le     C\rho+  CM_0^2t^{2s} +  C\rho^{-2}(1+M_0t^{s-1})(M_0^\theta t^{s\theta}+\|E\|_{H_0^1(D)}).
\ee
To finish the proof, we consider two cases.
\nl
{\bf Case 1:} $\|E\|_{H_0^1(D)}>0$.
First, we choose $t$ so that $M_0^\theta t^{s\theta}=\|E\|_{H_0^1(D)}$, i.e.   $ t:=\|E\|^{\frac{1}{s\theta} }_{H_0^1(D)}M_0^{-1/ s}$,   so that the two terms in the last
bracketed sum of \eref{finally1} are equal.       Since 
\be
\label{note1}
\|E\|_{H_0^1(D)}\le C ,
\ee
 and $M_0\ge C$ (because of    \eqref{e:bound_M}),  this choice of $t$ satisfies 
 \be
 \label{note2}
 1\leq CM_0t^{s-1}.
 \ee 
 Next,  we    choose $\rho$ such 
that $\rho^3=M_0t^{s-1}\|E\|_{H_0^1(D)}= M_0^{1/s}\|E\|^{\frac{s\theta +s-1}{s\theta}}_{H_0^1(D)}$.  This choice balances the first and last terms on the right side of \eref{finally1} and therefore gives

\be
\label{desired}
\|\delta\|^2_{L_2(D)}\le  CM_0^{\frac{1}{3s}} \|E\|^{\frac{s\theta+s-1}{3s\theta}}_{H_0^1(D)} +  C\|E\|^{\frac{2}{\theta}}_{H_0^1(D)}.
\ee
Since $\frac{s\theta+s-1}{3s}\le 2$,  the inequalities \eref{note1} and \eref{note2} show that  the first term in the sum on the right can be absorbed into the second, and the theorem follows.
\nl
{\bf Case 2:} $\|E\|_{H_0^1(D)}=0$.   For any  sufficiently small  $t>0$, we  choose $\rho$ such that $\rho^3=M_0^{1+\theta}t^{s\theta+s-1}$ so that the first and last terms in \eref{finally1} balance.  Then, 
  \eref{finally1} gives
 $$
 \|\delta\|^2_{L_2(D)}\le CM_0^{\frac{1+\theta}{3}}t^{\frac{s\theta+s-1}{3}}+CM_0^2t^{2s}.
 $$
Since by assumption,   $\theta>\frac{1-s}{s}$, we have $t^{s\theta+s-1}\rightarrow 0$ as $t\rightarrow 0$, and therefore \eref{3:t2} holds in this case as well.
 \hfill $\Box$
\nl

Note that the proof of the above  theorem relies on the fact that $(D,a_t,f)$ and $(D,b_t,f)$ both satisfy the positivity condition \lb{2} for a uniform constant $c$.
The proof  can be 
 easily modified to cover the case where  $(D,a_t,f)$ and $(D,b_t,f)$ satisfy the positivity condition \lb{\beta} with a uniform constant $c$  for  any given $0\leq \beta <2$.

%

%
%

  \begin{remark}
 \label{Remark:Falk}   
 As noted in the introduction,  a typical result based on least squares or variational techniques for finding the diffusion 
 coefficient $a$ is 
estimate \eref{est}. For clarity, we focus here on the results from \cite{KL,F}, 
where the  approximation $a_h\in A_h$
is computed solely based on the knowledge of $u^{ob}$. Therefore any two diffusion coefficients  $a$ and $b$ with the 
same observed $u^{ob}$ will have the same approximant $a_h$, generated by the above process.  If we  take $u^{ob}=u_a$ in \eref{est}, we obtain the bound  
\be
\|a-a_h\|_{L_2(D)}\le Ch^r.
\ee
On the other hand, we can view $u_a=u^{ob}$ as an observation of $u_b$ and in this case obtain from \eref{est}, the bound
\be
\label{rknow1}
\|b-a_h\|_{L_2(D)}\le C(h^r+h^{-2}\|u_a-u_b\|_{L_2(D)}).
\ee
Hence,
\be
\label{rknow2}
\|b-a\|_{L_2(D)}\le Ch^r+C(h^r+h^{-2}\|u_a-u_b\|_{L_2(D)}).
\ee
If  we chose $h$, such that  $h^r=h^{-2}\|u_a-u_b\|_{L_2(D)}$,  we obtain the estimate
\be
\label{know3}
\|b-a\|_{L_2(D)}\le C \|u_a-u_b\|^{\frac{r}{r+2}}_{L_2(D)}.
\ee
 Besides working with Neumann boundary conditions,
there are  two major distinctions between \eref{know3} and our results.   The first is the $L_2(D)$ norm that appears on the right side in place of our $H_0^1(D)$ norm.  Recall that we have already mentioned (see \eref{weaker}) how one can derive bounds of the form \eref{know3} from our results. The second distinction  is the much more demanding regularity assumption placed on $a,b$ as well as on $u_a,u_b$.  Namely,  \eref{know3} is proved in the above references under    the regularity requirements $a,b\in H^{r+1}(D)$ and     $u_a,u_b\in W^{r+3}(L_\infty(D))$
with $r\ge 1$.   Whereas, in our treatment, stability estimates are available solely under the much weaker stability assumption   $a,b\in H^s(D)$,  $s*<s\le 1$,
 where $s^*<1$. \end{remark}

\subsection{The gradient condition \grad{}}
\label{solution-smoothness}

The statement of Theorem \ref{5:istheorem1} relies on the assumption that the solutions $u_a$ and $u_b$ satisfy the gradient condition \grad{}. Finding sufficient conditions that ensure \grad{} is a well studied question in harmonic analysis and partial differential
equations.
We   recall,  two classes of diffusion coefficient for which such condition holds.

\subsubsection{VMO diffusion coefficients}

We start with the following result from \cite{AQ}.
\begin{result} 
If  $D$ is a   $C^{1}$ domain,  the diffusion coefficient $a$ is in ${\rm VMO}\cap \cA$,  and the right  side $f=\mathrm{div}(g)$, with $g\in L_p(D)$, 
then there exists a unique weak solution $u_{a}$ to \eref{i:equation} such that   $\nabla u_{a} \in L_p(D)$, $1< p<\infty$,  and
\be
\label{1:r1}
\|\nabla u_{a}\|_{L_p(D)}\le C\|g\|_{L_p(D)},
\ee
with $C$ depending only on $D, d, p, \lambda, \Lambda$ and the VMO modulus  of $a$.  
\end{result}

Recall that the VMO modulus $\nu(a,\cdot)$ of $a$ is defined by
$$
\nu(a,t):=\sup_{|Q|\le t} \frac{1}{|Q|}\int_Q|a-a_Q|, \quad a_Q:= \frac 1 {|Q|} \int_Q a,\quad t>0,
$$
where the supremum is taken over all cubes $Q$ with measure  at most $ t$. 
In order to show that $u_a$ satisfies the gradient condition \grad, we need to consider 
a subclass of diffusion coefficients $a$, for which the estimate \eref{1:r1} is uniform for all functions in this class. 
For this, we consider  a non-decreasing continuous function $\Phi(t)$, $t\ge 0$, with $\Phi(0)=0$,
and introduce the class  $\cA_{\Phi}$ defined as
\be
\label{1:VMObound}
\cA_{\Phi}:=\{a\in\cA\,:\,\nu(a,t) \le \Phi(t),\ t>0\}.
\ee
Likewise, for $s>0$, we define the class
\be
\label{1:VMObound1}
\cA_{s,\Phi}:=\cA_s\cap \cA_{\Phi}.
\ee

An examination of the proofs in \cite{AQ} and \cite{D} shows that for all $a\in\cA_{\Phi}$ the 
constant in \eref{1:r1} is uniformly bounded, with a bound, depending on $\Phi$, $D$, $d$, $\lambda$, $\Lambda$. Therefore, according to the estimate \eqref{1:r1},  for each $0<\theta<1$, the 
solution $u_a$ satisfies the gradient condition \grad{}  
 with $M$  only depending  on $\theta$, $D$, $d$, $\lambda$, $\Lambda$, $\Phi$, and $f$.
As a consequence, we deduce the following corollary of Theorem \ref{5:istheorem1}.
 \begin{corollary}
 \label{vmores}
Let $D$ be a $C^1$ domain, $f \in L_\infty(D)$ with $f\geq c_f>0$ and  $\Phi(t)$, $t\ge 0$, be a non-decreasing continuous function with $\Phi(0)=0$.
Furthermore, assume that $a,b\in \cA_{s,\Phi}$ for some $\frac 1 2<s\leq 1$. 
Then there exists a constant $C$ only depending on $D$, $d$, $\lambda$, $\Lambda$, $f$, and $\Phi$ such that
 \be
  \|a-b\|_{L_2(D)}\le C\sqrt{1+(\|a\|_{H^s(D)}+\|b\|_{H^s(D)})^{\frac {1} {3s}}} \|u_a-u_b\|_{H_0^1(D)}^{r}
 \ee
 for every   $r < \frac{2s-1}{6s}$. 
 \end{corollary}

 \hfill $\Box$

\subsubsection{General diffusion coefficients}

Again, we start with the following gradient estimate.
 \begin{result}[see \cite{M,BDN}]   
 If   $D$ is any Lipschitz domain, then  there is a value $P>2$, depending on $D$,  such that whenever  $a\in \cA$ and 
 $f\in W^{-1}(L_p(D))$, with $2\le p<P$, then 
%
$$
\|\nabla u_a\|_{L_p(D)}\le C\|f\|_{W^{-1}(L_p(D))},
$$
with $C$ depending only on $d,D,\lambda,\Lambda,p$.   
\end{result}
It follows from the above result that $u_a$ satisfies condition \grad{} for $0<\theta<\frac{P-2}{P}$,
 where $M$ depends on 
$d,D,\lambda,\Lambda$, and $f$.
Therefore,  Result 2 and Theorem \ref{5:istheorem1} lead to the following corollary.
\begin{corollary}
\label{cor4.4}
Let $D$ be a Lipschitz domain, $f\in L_\infty(D)$ with $f\geq c_f>0$ and let  $P>2$ be the constant   in   Result  2. 
Assume that $a,b\in \mathcal A_s$ with $\frac{P}{2(P-1)}<s\leq 1$.
Then, there exists a constant $C$ only depending on $D$, $d$, $s$, $\lambda$, $\Lambda$, and $f$ such that
 \be
  \|a-b\|_{L_2(D)}\le C\sqrt{1+(\|a\|_{H^s(D)}+\|b\|_{H^s(D)})^{\frac {1} {3s}}}\|u_a-u_b\|_{H_0^1(D)}^{r},
 \ee
 for every $r <\frac{1}{6}- \frac{P(1-s)}{6(P-2)s}$. 
 \end{corollary}
 


\section{Piecewise constant diffusion coefficients}
\label{s:pwc}

Piecewise constant diffusion coefficients are often used in numerical simulation.  
This case is not covered by the 
discussions  in the preceding sections because such diffusion coefficients  do not satisfy the regularity assumptions 
considered there.    In this section, we  
derive some elementary results for 
piecewise constant parameters $a$,  subordinate to a fixed partition.  
We assume for simplicity that the domain $D=(0,1)^d$ and $\cP_n$ is the 
partition of $D$ into $n^d$  disjoint cubes of side length $1/n$.   The  derivations that follow
can be generalized to other settings.  We denote by $\cA^n$  the set of all
diffusion coefficients $a$ defined on  $D$ that are piecewise constant functions subordinate to $\cP_n$.   
We continue to make the assumption that each $a\in\cA^n$ satisfies $\lambda\le a\le \Lambda$ for fixed  
$0<\lambda<\Lambda$, and therefore can be written as
\be
\label{5:an}
a:=\sum_{Q\in\cP_n} a_Q\chi_Q,
\ee
where  $a_Q\in [\lambda,\Lambda]$,  and $\chi_Q$ is the characteristic function of the cube $Q$.

\begin{lemma}
\label{5:regularitylemma}
Let $D=(0,1)^d$ and $f \in L_2(D)$. 
If the diffusion coefficient $a\in \cA^n$ is given by \eqref{5:an}, then  for each cube $Q\in\cP_n$,  
the solution $u_a$ to \eref{i:equation}  satisfies the equation
\be
\label{5:rep}
- a_Q \Delta u_a (x) =f(x), \quad  a.e.\   x\in Q.
\ee
\end{lemma}
\noindent
{\bf Proof:} Let $a\in\cA_n$ and  $Q\in\cP_n$. 
Following the proof of  the interior regularity theorem, see \cite{E}, one can show that 
$u_a\in W^2(L_2(\cO))$ on each open set  $\cO$ strictly contained in  $Q$.    If in  \eref{i:equation}, 
we take $v$ smooth and compactly supported
on $Q$ and integrate by parts, we find
\be
\label{5:ip}
-a_Q \int_Q \Delta u_a v=\int_Q fv.
\ee
It follows that $-a_Q \Delta u_a=f$ at every point $x$ in the interior of $Q$ which is a 
Lebesgue point of both $f$ and $\Delta u_a$.    In particular, this holds almost 
everywhere on $Q$. \hfill $\Box$

\begin{theorem}
\label{5:pwdtheorem}
Let $D=(0,1)^d$ and  $f\in L_2(D)$ with $f\ge c_f>0$ on $D$.   
Let $a,b\in\cA^n$ be diffusion coefficients and $u_a,u_b$ be the 
corresponding solutions to \eref{i:equation} on $D$.
Then for each $Q\in\cP_n$, we have
\be
\label{5:pwct1}
|a_Q-b_Q|\le C n^{\frac{d+2}{2}}\|\nabla u_a-\nabla  u_b\|_{L_2(Q)},
\ee
where $C$ depends only on $c_f$ and $\Lambda$.   Therefore, 
\be
\label{5:pwct2}
\|a-b\|_{L_2(D)}\le  Cn \|u_a-u_b\|_{H_0^1(D)}.
\ee

\end{theorem}

\noindent
{\bf Proof:}  From Lemma \ref{5:regularitylemma}, we know that for each $Q\in\cP_n$, we have
\be
\label{5:local}
a_Q-b_Q = \Delta (u_a-u_b)  \frac{a_Q b_Q}{f},\quad  {\rm a.e. \ on \ }Q.
\ee 
We now assume without loss of generality that $a_Q > b_Q$.
Therefore, we have that $\Delta(u_a - u_b )>0$ on $Q$ since $f>0$.
Recall that there exist functions 
$\varphi _Q\in C^\infty_c(Q)$ (for example the    standard mollifier  supported in $Q$),    
such that $\int_{Q} \varphi_Q = 1$ and
\be
\label{5:phi}
\|\nabla\varphi_Q\|_{L_2(Q)}\le C_0 n^{\frac{d+2}{2}},
\ee
 with $C_0$ an absolute constant.
Then multiplying \eref{5:local}  by  such a $\varphi_Q$ and integrating over $Q$ yields
$$
a_Q-b_Q  = \int_{Q} \Delta (u_a-u_b) \frac {a_Qb_Q} f \varphi_Q \leq \frac {a_Qb_Q} {c_f} \int_Q \Delta (u_a-u_b) \varphi_Q
= - \frac{a_Qb_Q} {c_f} \int_Q  \nabla (u_a-u_b) \nabla \varphi_Q,
$$
where we used integration by parts to get the last equality.
The boundedness of $a$ and $b$ yields
\begin{eqnarray}
\label{5:local1}
a_Q-b_Q  &\leq& C \| \nabla (u_a-u_b)\|_{L_2(Q)} \| \nabla  \varphi_Q  \|_{L_2(Q)}
\\ \nonumber
&\le&
Cn^{\frac{d+2}{2}} \| \nabla (u_a - u_b)\|_{L_2(Q)}.
\end{eqnarray}
This proves \eref{5:pwct1}.  To prove   \eref{5:pwct2}, we square \eref{5:pwct1} integrate over $Q$  to find
\be
\label{pcfind}
\int_Q|a-b|^2\le  Cn^{d+2} \| \nabla (u_a - u_b)\|^2_{L_2(Q)} n^{-d}= Cn^{2} \| \nabla (u_a - u_b)\|^2_{L_2(Q)}.
\ee
If we add these estimates up over all $Q\in\cP_n$ and take a square root, we arrive at \eref{5:pwct2}.
\hfill $\Box$


\section{The univariate case}\label{S:uni}

In the univariate case, several stability results, mainly for the  Neumann problem, are available, see for example, \cite{K}. 
Here, we will discuss the one dimensional Dirichlet problem  with diffusion coefficients  $a\in \cA$ 
and the domain $D=(0,1)$.  In this case, under certain assumptions on $f$, we will be able to improve the Lipschitz exponent in the inverse parameter estimate and also provide limits to how large this Lipschitz
exponent can be.

Notice that  in this case, one needs some assumptions on $f$ to guarantee that $a$ is uniquely 
determined from the solution $u_a$, as the following example, taken from \cite{K}, shows.
The function
 $$
   u(x) = \left\{ \begin{array}{ll} x, &  x \in [0, \frac 1 2], \\ 1-x, & x \in (\frac 1 2, 1],\end{array} \right.
 $$
is a solution on $D$ to the problem
$$
-(au')'=2\delta_{1/2}, \quad u(0)=u(1)=0,
$$
with diffusion coefficient $a\equiv 1$ or any  $a$ of the form
$$
   a= \left\{ \begin{array}{ll} q, &  \hbox{on}  \,\,[0, \frac 1 2], \\ 2-q, & \hbox{on}\,\, (\frac 1 2, 1],\end{array} \right.
 $$
where  $0<q<2$. Here $\delta_{1/2}$ is the delta distribution with weight $1$ at $1/2$. 

In going further, we consider the case $f=1$, noting that the  derivations below
can be generalized to other settings. We determine  the solution $u_a$    
and show that estimate \eref{3:it1} in Corollary \ref{3:cor} can be improved.
 We use the notation $A:=1/a$, $B:=1/b$, where $a,b\in\cA$.  Now,
 \eref{i:equation} becomes
\be
\label{weakone}
 \int_0^1 au_a'v'=\int_0^1  v,\quad  v\in H_0^1(0,1),
\ee
and one checks  that the solution to \eref{weakone}  is
\be
\label{onesol} 
u_a(x)=-\int_0^x A(t)(t-\gamma_a)\, dt,\quad \hbox{where}\quad \gamma_a:=\frac{\int_0^1A(t) t\, dt}{\int_0^1 A(t)\, dt} \in (0,1).
\ee
 This gives
\be
\label{v}
-A(x)(x-\gamma_a)=u_a'(x) .
\ee

\subsection{An upper bound}
\label{secub}
To bound  $\|a-b\|_{L_2(0,1)}$ in terms of $\|u_a'-u_b'\|_{L_2(0,1)}$, it is sufficient to   bound   $\|A-B\|_{L_2(0,1)}$.
  Let us set  $\eta:=\gamma_a - \gamma_b$, 
and   $E'(x):=u_a'(x)-u_b'(x)$. Without loss of generality, we may assume that $\eta\geq 0$, since  
otherwise we can reverse the roles of $a$ and $b$.  
The following lemma gives an estimate for $\eta$.
\begin{lemma}
\label{eta} We have 
\be
\label{lower2}
\eta\le  c_0\|E'\|^{2/3}_{L_2(0,1)}, 
\ee
where the constant $c_0$ depends only on $\lambda$ and $\Lambda$.
\end{lemma}

\noindent
{\bf Proof:}
The estimate obviously holds if  $\eta=0$, so we assume  that $\eta>0$. We consider an interval $I$ of 
length $2c\eta$ centered at $\gamma_a$ with 
$c:=\frac{\lambda}{2(\lambda+\Lambda)}<1/2$.  We have for  $x\in I\cap (0,1)$
\begin{eqnarray}
\nonumber
|u_a'(x)-u_b'(x)|&=& |(x- \gamma_b)B(x) -(x-\gamma_a )A(x)|\ge (1-c)B(x)\eta  -c A(x)\eta
\\ \nonumber
&\ge&\left(\frac{1-c}{\Lambda}-\frac{c}{\lambda}\right)\eta= \frac{\eta}{2\Lambda}.
\nonumber
\end{eqnarray}
Squaring this estimate  and integrating over $I\cap (0,1)$ gives
%
$$
\frac{\eta^2}{4\Lambda^2} |I\cap (0,1)|\le \|u_a'-u_b'\|_{L_2(0,1)}^2=\|E'\|_{L_2(0,1)}^2,
$$
and since $|I\cap (0,1)|\geq c\eta$,  the proof is completed.
\hfill $\Box$

The following lemma gives an upper bound for the norm $\|A-B\|_{L_2(0,1)}$.
\begin{lemma}
\label{auxlemma}
For every  $\rho>0$, we have
\begin{equation}
\label{tt}
\|A-B\|^{2}  _{L_2(0,1)}\le  \frac{C}{\rho^2}\|E'\|^{4/3}_{L_2(0,1)}+8\lambda^{-2}\rho,
\end{equation}
where $C$ depends only on $\lambda$ and $\Lambda$.
In particular, if $\|E'\|_{L_2(0,1)}=0$, then $A=B$ a.e in $(0,1)$.
\end{lemma}

\noindent
 {\bf Proof:} First, let us observe that
 \be
\label{start1}
(A(x)-B(x))(x-\gamma_a)=A(x)(x-\gamma_a)-B(x)(x-\gamma_b)
+B(x)(\gamma_a-\gamma_b)=-E'(x)+B(x)(\gamma_a-\gamma_b).
\ee
We now consider an interval $J$ of length $2\rho $ centered at $ \gamma_a $.  Then,  using \eref{start1} on $J^c$, where $J^c$ 
is the complement  of $J$ in $(0,1)$
(which might be empty), we have
$$
\rho|(A(x)-B(x)|\leq   |E'(x)|+\lambda^{-1}\eta,\quad x\in J^c,
$$
and therefore
$$
\rho^2|(A(x)-B(x)|^2\leq2|E'(x)|^2+2\lambda^{-2}\eta^2, \quad x\in J^c.
$$
We integrate the latter inequality over $J^c$ to obtain
\be
\label{comp}
 \|A-B\|^{2}  _{L_2(J^c)}\le \frac{2}{\rho^2}\|E'\|^{2}_{L_2(0,1)}   +  \frac{2\lambda^{-2}}{\rho^2}\eta^2.
\ee
Meanwhile, for $x\in J \cap (0,1)$, we have  $|A(x)-B(x)|\le 2\lambda^{-1}$ and therefore
\be
\label{onJ}
\|A-B\|^2_{L_2(J \cap (0,1))}\le 8\lambda^{-2}\rho.
\ee
Combining this with \eref{comp}, we obtain
\begin{eqnarray}
\label{final}
\|A-B\|^{2}_{L_2(0,1)}&\le & \frac{2}{\rho^2}\|E'\|^{2}_{L_2(0,1)} + \frac{2\lambda^{-2}}{\rho^2}\eta^2 +8\lambda^{-2}\rho
\nonumber\\
&\leq& \frac{2}{\rho^2}\|E'\|^{2}_{L_2(0,1)} + \frac{2c_0^2}{\rho^2\lambda^{2}}\|E'\|^{4/3}_{L_2(0,1)}+8\lambda^{-2}\rho,
\nonumber
\end{eqnarray}
where we have used Lemma \ref{eta}. Since  
$|u_a'(x)-u_b'(x)|= |(x- \gamma_b)B(x) -(x-\gamma_a )A(x)|\leq 2\lambda^{-1}$, we have that 
$\|E'\|_{L_2(0,1)}\leq  2\lambda^{-1}$, and the first term of the above inequality is absorbed by the second term. Hence, we get
$$
\|A-B\|^{2}_{L_2(0,1)}\leq \frac{C}{\rho^2}\|E'\|^{4/3}_{L_2(0,1)}+8\lambda^{-2}\rho,
$$
where $C$ depends only on $\lambda$ and $\Lambda$. This  proves the first part of the lemma. When $\|E'\|_{L_2(0,1)} =0$, 
$$
\|A-B\|^{2}_{L_2(0,1)}\le 8\lambda^{-2}\rho,
$$
for all $\rho>0$ and so $A=B$ a.e. in $(0,1)$.
\hfill $\Box$

We can now prove the following stability estimate in the one dimensional case.
\begin{theorem}
\label{thd1}
For any $a,b\in\cA$, the solutions $u_a,u_b$ to {\rm\eref{i:equation}} with $f=1$
satisfy the estimate
\be
\label{maind1}
\|a-b\|_{L_2(0,1)}\leq C\|u_a-u_b\|^{2/9}_{H^1_0(0,1)},
\ee
where $C$ depends only on $\lambda$ and $\Lambda$. 
In particular, if $u_a=u_b$ on $(0,1)$, then
$a=b$  a.e in $(0,1)$.
\end{theorem}

\noindent
{\bf Proof:}
If $\|u_a-u_b\|_{H^1_0(0,1)}=0$, it follows from Lemma \ref{auxlemma} that $a=b$, a.e. on $(0,1)$, and therefore \eref{maind1} holds.
When $\|E'\|_{L_2(0,1)}=\|u_a' - u_b'\|_{L_2(0,1)}>0$, we choose $\rho=\|E'\|^{4/9}_{L_2(0,1)}$ in Lemma \ref{auxlemma} to derive the desired estimate.
\hfill $\Box$

\subsection{A lower bound}
In this section, we show that the exponent in estimates of the form \eref{maind1}  cannot be greater than $1/3$. 

\begin{theorem}
\label{ctex}
Consider equation {\rm \eref{i:equation}} with domain $D=(0,1)$ and right  side $f=1$. There are diffusion coefficients $a,b\in \cA$, such that 
the corresponding solutions $u_a, u_b$, satisfy the inequality
\be
\label{final2}
\|a-b\|_{L_2(D)}\ge c\|u_a-u_b\|_{H_0^1(D)}^{1/3},
\ee
where $c$ is a constant, depending only  on $\lambda$ and $\Lambda$.
\end{theorem}
\noindent
{\bf Proof:} We define the following diffusion coefficients 
\begin{eqnarray*}
\frac{1}{a(x)}:=A(x)&=& \begin{cases}1, &\mbox{ for  } 0<x\le \alpha,\\
 2,& \mbox{ for  } \alpha<x<1,
\end{cases}\\ \\
\frac{1}{b(x)}:=B(x)&=& \begin{cases}1, &\mbox{ for  } 0<x\le \beta,\\
 2,& \mbox{ for  } \beta<x<1,
\end{cases}
\end{eqnarray*}
where $\alpha,\beta\in(0,1)$, and compute 
\be
\label{AB} 
\|A-B\|_{L_2(0,1)}= |\alpha-\beta|^{1/2}.
\ee
Let $g(t):= \frac{1-t^2/2}{2-t}$.  Then, a simple calculation gives
\be
\label{gamma}
\gamma_a= g(\alpha), \qquad \gamma_b = g(\beta),
\ee
where $\gamma_a$ and $\gamma_b$ are defined by \eref{onesol}.  We  denote by  $\alpha_0$  the point where $g$ achieves its minimum in $(0,1)$. 
Then, we have 
\be
\label{gder}  g'(\alpha_0)= 1-2\alpha_0+\alpha_0^2/2=0 \quad {\rm and} \quad  \alpha_0=2\sqrt{2}-2.  
\ee
  We   fix $\alpha$ as $\alpha_0$.  Since  $g(\alpha_0)=\alpha_0$, we have $\gamma_a=\alpha_0$. 
  
We now   bound $\eta:=\gamma_a-\gamma_b=\alpha_0-\gamma_b$ from above. 
In fact,  using \eref{AB} and \eref{gder}, we have
\be
\label{pop}
|\eta|= g(\beta)-g(\alpha_0)=\frac{(\alpha_0-\beta)^2}{2(2-\beta)}<\frac{1}{2}(\alpha_0-\beta)^2=\frac{1}{2}\|A-B\|^4_{L_2(0,1)}.
\ee
 Recall that
\be
\label{know1}
E'(x)=-(A(x)-B(x))(x-\gamma_a)+B(x)(\gamma_a-\gamma_b) =-(A(x)-B(x))(x-\alpha_0)+B(x)\eta.
\ee
Therefore, using \eref{AB} and \eref{pop} ,  we have
\begin{eqnarray}
\label{know3one}
\|E'\|^2_{L_2(0,1)} &\le& 2\int_0^1(A(x)-B(x))^2(x-\alpha_0)^2\, dx+2\eta^2 \int_0^1B^2(x)\, dx\nonumber\\
&\le&2\left |\int_{\alpha_0}^\beta(x-\alpha_0)^2\,dx\right |+8\eta^2= \frac{2}{3}|\beta-\alpha_0|^3+8\eta^2=
 \frac{2}{3}\|A-B\|^6_{L_2(0,1)}+8\eta^2
 \nonumber \\
 &\leq& \frac{2}{3}\|A-B\|^6_{L_2(0,1)}+2\|A-B\|^8_{L_2(0,1)}\leq C\|A-B\|^6_{L_2(0,1)},
\end{eqnarray}
where  $C$ depends only on $\lambda,\Lambda$.  This completes the proof.
\hfill $\Box$


\vskip 1in

 \noindent
 Andrea Bonito\\
Department of Mathematics, Texas A\&M University,
College Station, TX 77840, USA\\
  bonito@math.tamu.edu
\vskip .1in

 \noindent
 Albert Cohen\\
  Laboratoire Jacques-Louis Lions,
    UPMC Univ Paris 06, UMR 7598, F-75005, Paris, France\\
cohen@ann.jussieu.fr

\vskip .1in

\noindent
Ronald DeVore\\
Department of Mathematics, Texas A\&M University,
College Station, TX 77840, USA\\
  rdevore@math.tamu.edu

 \vskip .1in
\noindent
Guergana Petrova\\
Department of Mathematics, Texas A\&M University,
College Station, TX 77840, USA\\
gpetrova @math.tamu.edu

 \vskip .1in
\noindent
Gerrit Welper\\
Department of Mathematics, Texas A\&M University,
College Station, TX 77840, USA\\
  gwelper@math.tamu.edu
\\

\end{document}